\documentclass[11pt]{article}

\usepackage[T2A]{fontenc}
\usepackage[cp1251]{inputenc}
\usepackage{amsfonts}
\usepackage{eufrak}
\usepackage{amssymb}
\usepackage{amsmath}
\usepackage{cite}



\begin{document}

\newtheorem{theorem}{Theorem}
\newtheorem{lemma}{Lemma}
\newtheorem{corollary}{Corollary}
\newtheorem{definition}{Опеределение}
\newtheorem{proposition}{Proposition}
\newtheorem{remark}{Remark}

\bigskip

\centerline{\textbf{The Heyde theorem on a group $\mathbb{R}^n\times
D$, where $D$ is a discrete Abelian group}}

\bigskip

\centerline{\textbf{Margaryta Myronyuk}}

\bigskip

\centerline{\textit{B. Verkin Institute for Low Temperature Physics
and Engineering }}

\centerline{\textit{of the National Academy of Sciences of Ukraine,
}}

\centerline{\textit{47 Nauky Ave, Kharkiv, 61103, Ukraine}}

\centerline{{Email: myronyuk@ilt.kharkov.ua}}

\begin{abstract}
Heyde proved that a Gaussian distribution on the real line is
characterized by the symmetry of the conditional distribution of one
linear statistic given another. The present article is devoted to a
group analogue of the Heyde theorem. We describe distributions of
independent random variables $\xi_1$, $\xi_2$ with values in a group
$X=\mathbb{R}^n\times D$, where $D$ is a discrete Abelian group,
which are characterized by the symmetry of the conditional
distribution of the linear statistic $L_2 = \xi_1 + \delta\xi_2$
given $L_1 = \xi_1 + \xi_2$, where $\delta$ is a topological
automorphism of $X$ such that ${Ker}(I+\delta)=\{0\}$.

\end{abstract}



\bigskip

\section{Introduction}

Characterization problems in mathematical statistics are statements
in which the description of possible distributions of random
variables follows from properties of some functions in these
variables. A large number of studies have been devoted to theorems
characterizing the Gaussian distribution on the real line. In
particular, in 1970 C.C.Heyde proved the following theorem.

\medskip

\textbf{The Heyde theorem} (\cite{He}, \cite[\S\,13.4.1]{KaLiRa}).
\textit{Let $\xi_j, \ j=1, 2, ..., n, \ n \ge 2,$ be independent
random variables. Consider linear statistics $L_1=\alpha_1\xi_1+
\cdots +\alpha_n\xi_n$ and $L_2=\beta_1\xi_1 + \cdots +
\beta_n\xi_n$, where the coefficients $\alpha_j, \ \beta_j $ are
nonzero real numbers such that $\beta_i\alpha_i^{-1} \pm
\beta_j\alpha_j^{-1} \ne 0$ for all $i \ne j$. If the conditional
distribution of $L_2$ given $L_1$ is symmetric then all random
variables $\xi_j$ are Gaussian.}

\medskip

The theory of characterization theorems is being actively developed
in a situation when random variables take values in a locally
compact Abelian group (see e.g. \cite{Fe5a}). Group analogues of the
Heyde theorem were studied in the case when independent random
variables take values in locally compact Abelian groups and
coefficients of the linear statistics are either topological
automorphisms of the group (\cite{Fe2, Fe4, Fe3, Fe6, F, Fe7, My2,
My1, MyF}) or integers (\cite{My2020}). Basically, the following
problem was studied.

Let $X$ be a second countable locally compact Abelian group. Denote
by ${\rm Aut}(X)$ the group of topological automorphisms of the
group $X$ and denote by $I$ the identity automorphism. Let $\xi_1$
and $\xi_2$ be independent random variables with values in $X$ and
distributions $\mu_1$ and $\mu_2$. Let $\delta \in {\rm Aut}(X)$ and
\begin{equation}\label{i1}
I\pm\delta\in{\rm Aut}(X).
\end{equation}
What can be said about the distributions $\mu_1$ and $\mu_2$, if the
conditional distribution of the linear statistic $L_2 = \xi_1 +
\delta\xi_2$ given $L_1 = \xi_1 + \xi_2$ is symmetric?

This problem was studied on discrete groups, the real
$m$-dimensional space, products of these groups, $\textbf{a}$-adic
solenoids. The proofs of the obtained results essentially used
condition (\ref{i1}).

It turns out (see \cite{FeTVP1, FeTVP2, Fe7}) that in some cases
theorems, which were proved under condition (\ref{i1}), remain true
if condition (\ref{i1}) is replaced by the weaker condition

\begin{equation}\label{i2}
{\rm Ker}(I+\delta)=\{0\}.
\end{equation}
In this article we generalize results of \cite{FeTVP2} and describe
distributions on a group of the form

\begin{equation}\label{i3}
X=\mathbb{R}^n\times D,
\end{equation}
where $D$ is a discrete Abelian group, which are characterized by
the symmetry of the conditional distribution of one linear statistic
given another. Note that our result also generalizes the result of
the article \cite{My2}, where the group analogue of the Heyde
theorem was proved on the group of form (\ref{i3}), where $n=1$ and
condition (\ref{i1}) holds.

\section{Notation and definitions}

In the article we use standard results on the structure theory of
locally compact Abelian groups and abstract harmonic analysis (see
e.g. \cite{Hewitt-Ross}). Let $X$ be a second countable locally
compact  Abelian group, $Y=X^\ast$ be its character group, and
$(x,y)$ be the value of a character $y \in Y$ at an element $x \in
X$. Let $K$ be a subgroup of $Y$. Denote by $A(X,K)=\{x \in X:
(x,y)=1 \ \ \forall \ y \in K\}$ the annihilator of $K$. Denote by
$b_X$ the subgroup of all compact elements of $X$, and denote by
$c_X$ the connected component of zero of $X$. If $\delta$ is a
continuous endomorphism of the group $X$ then the adjoint
endomorphism $\tilde\delta$ of the group $Y$ is defined by the
formula $(x, \tilde\delta y) = (\delta x, y)$ for all $x \in X$, $y
\in Y$. Let $p$ be a prime number. The $p$-component of an Abelian
group is a subgroup consisting of elements whose order is a power of
$p$. Denote by $X_{p}$ the $p$-component of the group $X$. A torsion
group $X$ is called $p$-primary, if $X=X_{p}$. For each integer $n$,
$n \ne 0,$ let $f_n : X \mapsto X$ be the endomorphism $f_n x=nx.$
Set $X^{(n)} = f_n(X)$, $X_{(n)}=Ker f_n$. A subgroup $G$ of a group
$X$ is said to be characteristic if $G$ is invariant under each
topological automorphism of a group $X$. The subgroups $b_X$, $c_X$,
$X^{(n)}$, $X_{(n)}$ are characteristic.

Let $f(y)$ be a function on $Y$, and $h\in Y.$ Denote by $\Delta_h$
the finite difference operator
    $$\Delta_h f(y)=f(y+h)-f(y), \quad y \in Y.$$
A function $f(y)$ on $Y$ is called a polynomial if
    $$\Delta_{h}^{n+1}f(y)=0$$
for some $n$ and for all $y,h \in Y$.

Let ${M^1}(X)$ be the convolution semigroup of probability
distributions on $X$. Let $\hat \mu(y) = \int_X (x, y) d\mu(x)$ be
the characteristic function of a distribution $\mu \in {M^1}(X)$,
and $\sigma(\mu)$ be the support of $\mu$. Define $\bar \mu \in {\rm
M}^1(X)$ by the formula $\bar \mu(B) = \mu(-B)$ for all Borel sets
$B$ in $X$. Then $\hat{\bar{\mu}}(y)=\overline{\hat\mu(y)}$. Denote
by $m_K$ the Haar distribution of a compact subgroup $K$ of a group
$X$ and denote by $E_x$ the degenerate distribution concentrated at
the point $x\in X$. Note that the characteristic function of the
distribution $m_K$ has the form
\begin{equation}\label{11a}
\hat m_K(y)=
\begin{cases}
1, & \text{\ если\ }\   y\in A(Y, K),
\\  0, & \text{\ если\ }\ y\not\in
A(Y, K).
\end{cases}
\end{equation}
Denote by $\Gamma({\mathbb R}^n )$ the set of Gaussian distributions
on ${\mathbb R}^n $.

\section{Main theorem}

The main result of the article is the following theorem.

\begin{theorem}\label{m}
    Let $X=\mathbb{R}^n\times D$, where $D$ is a countable discrete Abelian group.
    Let $\delta\in Aut(X)$ such that condition $(\ref{i2})$ is fulfilled.
    Let $\xi_1$ and $\xi_2$ be independent random variables with values in $X$ and distributions $\mu_1$ and $\mu_2$.
    If the conditional distribution of the linear statistic $L_2 = \xi_1 + \delta\xi_2$ given $L_1 = \xi_1 + \xi_2$ is symmetric, then
    $\mu_j=\gamma_j*\rho_j*m_F*E_{g_j}$, where $\gamma_j \in \Gamma({\mathbb R}^n)$, $\sigma(\rho_j)\subset D_2$,
    $F$ is a finite subgroup of $D$ without elements of order 2, $g_j\in D$, $j=1, 2$. Moreover, $\delta(F)=F$.
\end{theorem}

As appears from the paper \cite{MyF}, on a 2-primary finite group,
even under condition (\ref{i1}), one can hardly expect to obtain a
reasonable description of distributions characterized by the
symmetry of the conditional distribution of one linear statistic
given another.

We need some lemmas to prove Theorem \ref{m}.

\begin{lemma}\label{l1}
    (\cite[Lemma 16.1]{Fe5a})
Let $X$ be a second countable locally compact Abelian group,
$Y=X^*$. Let $\delta\in Aut(X)$. Let $\xi_1$ and $\xi_2$ be
independent random variables with values in $X$ and distributions
$\mu_1$ and $\mu_2$. The conditional distribution of the linear
statistic $L_2 = \xi_1 + \delta\xi_2$ given $L_1 = \xi_1 + \xi_2$ is
symmetric if and only if the characteristic functions $\hat\mu_j(y)$
satisfy the equation
\begin{equation}\label{l2}
\hat\mu_1(u+v )\hat\mu_2(u+\varepsilon v )= \hat\mu_1(u-v
)\hat\mu_2(u-\varepsilon v), \quad u, v \in Y,
\end{equation}
where $\varepsilon=\tilde{\delta}$.
\end{lemma}

\begin{lemma}\label{pol} (\cite{Fe-MarcLuk}) Let $Y$ be a locally compact Abelian group, $f(y)$ be a continuous polynomial on $Y$.
Then $f(y)=const$ for $y\in b_Y.$
\end{lemma}

For convenience, we formulate the following well-known statements.

\begin{lemma}\label{supp1} Let $X$ be a second countable locally compact Abelian group, $Y=X^*$.
Let  $\mu\in{\rm M}^1(X)$. Then the set $E=\{y\in Y:\
\hat\mu(y)=1\}$ is a closed subgroup of $Y$ and $\sigma(\mu)\subset
A(X,E)$.
\end{lemma}

\begin{lemma}\label{supp} (\cite[Lemma 2.13]{Fe5a}) Let $X$ be a topological group, $G$ be a Borel subgroup of
$X$, $\mu\in {\rm M}^1(G)$, $\mu=\mu_1*\mu_2$, where $\mu_j\in {\rm
M}^1(X)$. Then the distributions $\mu_j$ can be replaced by their
shifts $\mu'_j$ in such a manner that $\mu=\mu'_1*\mu'_2$ and
$\mu'_j\in {\rm M}^1(G)$.
\end{lemma}

The following lemma is crucial for the proof of Theorem \ref{m}.

\begin{lemma}\label{t1}
    Let $X=\mathbb{R}^n\times G$, where $G$ is a countable discrete 2-primary Abelian group.
    Let $\delta\in Aut(X)$ such that condition $(\ref{i2})$ is fulfilled.
    Let $\xi_1$ and $\xi_2$ be independent random variables with values in $X$ and distributions $\mu_1$ and
    $\mu_2$. If the conditional distribution of the linear statistic $L_2 = \xi_1 + \delta\xi_2$ given $L_1 = \xi_1 + \xi_2$ is symmetric,
    then
    $\mu_j=\gamma_j*\rho_j$, where $\gamma_j \in \Gamma({\mathbb R}^n)$, $\sigma(\rho_j)\subset G$, $j=1, 2$.
\end{lemma}

\textbf{Proof.} The group $Y=X^*$ is topologically isomorphic to the
group ${\mathbb R}^n\times H$, where $H=G^*$. To avoid introducing
new notation we will suppose that $Y={\mathbb R}^n\times H$. Denote
by $x=(t, g)$, where $t\in {\mathbb R}^n$, $g\in G$, elements of the
group $X$. Denote by $y=(s, h)$, where $s\in {\mathbb R}^n$, $h\in
H$, elements of the group $Y$. Put $\varepsilon=\tilde{\delta}$.

\textbf{1.} First we shall show that the restriction of the
endomorphism $I-\delta$ on $G$ has a zero kernel. If $x\in G_{(2)}$,
$x\neq 0$, then $(I-\delta)x=(I+\delta)x$. It follows from the
condition (\ref{i2}) that $x\not \in Ker(I-\delta)$. Suppose that
there exists an element $x_0\in G\setminus G_{(2)}$ such that
$(I-\delta)x_0=0$. Since the subgroup $G$ is 2-primary, the element
$x_0$ has an order $2^k$ for some $k>1$, i.e. $2^k x_0=0$ and
$2^{k-1} x_0\neq 0$. Since $(I-\delta)x_0=0$, we have
$(I-\delta)2^{k-1} x_0=0$. Since $2^{k-1} x_0\in G_{(2)}\setminus
\{0\}$, as has been shown above, $2^{k-1} x_0\not \in
Ker(I-\delta)$. We obtain the contradiction. Thus, we get
\begin{equation}\label{t3}
    (I-\delta)x\neq 0, \quad x\in G\setminus \{0\}.
\end{equation}

Since the subgroups $G$ and $H$ are characteristic, conditions
(\ref{i2}), (\ref{t3}) and the compactness of $H$ implies that

\begin{equation}\label{t4}
    (I\pm\varepsilon)H=H.
\end{equation}

\textbf{2.} We reduce the proof of the lemma to the case when the
subgroup $G$ is bounded. It follows from Lemma \ref{l1} that the
characteristic functions $\hat\mu_j(y)$ satisfy equation (\ref{l2}).
Put $\nu_j = \mu_j * \bar \mu_j$. Then $\widehat \nu_j(y) =
|\widehat \mu_j(y)|^2 \ge 0$, $y \in Y$. It is obvious that the
characteristic functions $\widehat \nu_j(y)$ satisfy equation
(\ref{l2}) too. Since the subgroup $H$ is characteristic, we can
consider the restriction of equation (\ref{l2}) on the subgroup $H$.

Let $U$ be a neighborhood of zero of $H$ such that all
characteristic functions $\nu_j(y)> 0$ for $y\in U$. Put
$\varphi_j(y) = - \ln \widehat\nu_j(y), \ y \in U$.

Let $V$ be a neighborhood of zero of $H$ such that

$$
\sum_{j=1}^8 \lambda_j (V) \subset U
$$
for all endomorphisms $\lambda_j \in \{I, \varepsilon \}$.

Since $G$ is a discrete torsion group, the group $H$ is compact
totally disconnected. Then each  neighborhood of zero of $H$
contains an open compact subgroup (\cite[7.7]{Hewitt-Ross}). Let
$W\subset V$, where $W$ is an open subgroup of $H$. We prove that
the functions $\varphi_j(y)$ are polynomials on a subgroup.

It follows from (\ref{l2}) that the functions $\varphi_j(y)$ satisfy
equation

\begin{equation}\label{l2-1}
\varphi_1(u+v) + \varphi_2(u+\varepsilon v)- \varphi_1(u- v) -
\varphi_2(u-\varepsilon v)=0, \quad u, \ v \in W.
\end{equation}

We use the finite differences method. Let $k_1$ be an arbitrary
element of $W$. Substitute in (\ref{l2-1}) $u+\varepsilon k_1$ for
$u$ and $v+k_1$ for $v$. Subtracting equation (\ref{l2-1}) from the
resulting equation we obtain

\begin{equation}\label{l2-2}
\Delta_{l_{11}}\varphi_1(u+ v) +
\Delta_{l_{12}}\varphi_2(u+\varepsilon v) -
\Delta_{l_{13}}\varphi_1(u- v)=0, \quad u, v \in W,
\end{equation}

\noindent where $l_{11}= (I+\varepsilon)k_1$, $l_{12}=2 \varepsilon
k_1$, $l_{13}= (\varepsilon-I)k_1$. Let $k_2$ be an arbitrary
element of $W$. Substitute in (\ref{l2-2}) $u+k_2$ for $u$ and
$v+k_2$ for $v$. Subtracting equation (\ref{l2-2}) from the
resulting equation we obtain

\begin{equation}\label{l2-3}
\Delta_{l_{21}}\Delta_{l_{11}}\varphi_1(u+ v) +
\Delta_{l_{22}}\Delta_{l_{12}}\varphi_2(u+\varepsilon v) = 0,
 \quad u, v \in W,
\end{equation}

\noindent where $l_{21}=2  k_2$, $l_{22}= (I+\varepsilon)k_2$. Let
$k_3$ be an arbitrary element of $W$. Substitute in (\ref{l2-3})
$u-\varepsilon k_3$ for $u$ and $v+k_3$ for $v$. Subtracting
equation (\ref{l2-3}) from the resulting equation we obtain

\begin{equation}\label{l2-4}
\Delta_{l_{31}}\Delta_{l_{21}}\Delta_{l_{11}}\varphi_1(u+ v) = 0,
\quad u, v \in W,
\end{equation}

\noindent where $l_{31}= (I-\varepsilon)k_3$. Putting $v=0$ in
(\ref{l2-4}), we get

\begin{equation}\label{l2-5}
\Delta_{l_{31}}\Delta_{l_{21}}\Delta_{l_{11}}\varphi_1(u) = 0, \quad
u \in W.
\end{equation}

Since elements $k_j$ are arbitrary, it follows from (\ref{t4}),
(\ref{l2-5}) and the expressions of $l_{11}, \ l_{21}, l_{31}$ that
the function $\varphi_1(y)$ satisfy the equation

\begin{equation}\label{l2-6}
\Delta_h^3\varphi_1(y)=0, \quad h, y \in B,
\end{equation}
\noindent where the subgroup $B=(I+\varepsilon)W \cap
(I-\varepsilon)W \cap W^{(2)}$. Thus the function $\varphi_1(y)$ is
a polynomial on the subgroup $B$.

Similarly, we get that the function $\varphi_2(y)$ satisfy equation
(\ref{l2-6}) too.

Since the subgroup $W$ is open, it is closed and therefore compact.
Hence the subgroup $B$ is also compact. It follows from this, Lemma
\ref{pol} and the condition $\widehat\nu_j(0)=1$ that
$\varphi_j(y)=0$ for $y\in B$. Hence, $\widehat\nu_j(y)=1$ for $y\in
B$. Lemma \ref{supp1} implies that $\sigma(\nu_j)\subset A(X,B)$.
Put $D=A(X,B)$. Note that the subgroup $D$ is generated by the
annihilators $A(X,(I+\varepsilon)W \cap (I-\varepsilon)W)$ and
$A(X,W^{(2)})$.

Put $W'= (I+\varepsilon)W \cap (I-\varepsilon)W$. It follows from
condition (\ref{t4}) that the endomorphisms $I+\varepsilon,
I-\varepsilon$ are open. Then the subgroup $W'$ is also open in $H$.
Note that $(Y/W')^*\approx A(X,W')$. Since the factor-group $H/W'$
is finite, the annihilator $A(X,W')=A(\mathbb{R}^n\times
G,W')=\mathbb{R}^n\times A(G,W') = \mathbb{R}^n\times F'$, where
$F'$ is a finite subgroup of $G$.

We have $$A(X,W^{(2)})=\{ x\in X: (x,2w)=1 \ \forall w\in W \}=\{
x\in X: (2x,w)=1 \ \forall w\in W \}= $$$$= \{ x\in X: 2x\in A(X,W)
\}= f_2^{-1}\left(A(X,W)\right).$$

Note that $Ker f_2=X_{(2)}=G_{(2)}$. Since $(Y/W)^*\approx A(X,W)$
and $W$ is an open subgroup of $H$, it follows from the finiteness
of the factor group $H/W$ that the annihilator
$A(X,W)=\mathbb{R}^n\times F''$, where $F''$ is a finite subgroup of
$G$. We get that the subgroup $A(X,W^{(2)})$ is generated by
$G_{(2)}$ and the subgroup $\mathbb{R}^n\times F''$.

Thus, the subgroup $D$ is generated by the subgroups
$\mathbb{R}^n\times F'$, $\mathbb{R}^n\times F''$ and $G_{(2)}$.
Hence $D\subset \mathbb{R}^n\times G_{(k)}$ for some $k$. So we get
that $\widehat\nu_j(y)=1$ for $y\in B$. Lemma 3 implies that
$\sigma(\nu_j)\subset D\subset \mathbb{R}^n\times G_{(k)}$ for some
$k$.

Lemma \ref{supp} implies that the distributions $\mu_j$ are
concentrated on the sets $x_j+(\mathbb{R}^n\times G_{(k)})$ for some
$x_j\in X$. Since $x_j=(t_j,g_j)$, where $2^{l_j}g_j=0$ for some
nonnegative integers $l_j$, there exists a nonnegative integer $l$
such that the supports $\sigma(\mu_j)\subset \mathbb{R}^n\times
G_{(l)}$. Since the subgroup $\mathbb{R}^n\times G_{(l)}$ is
characteristic, we can prove the lemma in the case when
$X=\mathbb{R}^n\times G_{(l)}$.

\textbf{3.} So, let $X=\mathbb{R}^n\times G$, where $G$ ia a
discrete bounded 2-primary group. Then $Y=\mathbb{R}^n\times H$,
where $H$ is a compact bounded 2-primary group. It is obvious that
conditions (\ref{i2}), (\ref{t3}) and (\ref{t4}) are fulfilled.
Further reasoning is similar to the reasoning of the article
\cite{MyF}.

Write equation (\ref{l2}) in the form

\begin{equation}\label{f1}
\hat\mu_1(s+s', h+h')\hat\mu_2(s+\varepsilon s', h+\varepsilon h')
=$$ $$ = \hat\mu_1(s-s', h-h')\hat\mu_2(s-\varepsilon s',
h-\varepsilon h'), \quad\quad (s,h), \ (s',h') \in Y.
\end{equation}

Put $h=h'=0$ in (\ref{f1}). We get

$$ \hat\mu_1(s+s',0)\hat\mu_2(s+\varepsilon s',0) = \hat\mu_1(s-s',0)\hat\mu_2(s-\varepsilon
s',0), \quad\quad s,\ s' \in \mathbb{R}^n.  $$

\noindent It was proved in \cite{FeTVP1} that all solutions of this
equation are the characteristic functions of Gaussian distributions,
i.e.

\begin{equation}\label{f2}
 \hat\mu_1(s, 0) = \exp\{-\langle A_1 s,s\rangle + i \langle t_1, s\rangle\}, \quad \hat\mu_2(s, 0) =
\exp\{-\langle A_2 s,s\rangle + i \langle t_2, s\rangle\},\quad s
\in \mathbb{R}^n,
\end{equation}
where $A_j \ge 0$ are positive semidefinite matrices, $t_j \in
\mathbb{R}^n$, $j = 1,2.$

We will prove by induction by $k$, where $2^k$ is the order of an
element $h$, that
\begin{equation}\label{f3}
\hat\mu_1(s,h)=\phi_1(s)\psi_1(h), \quad
\hat\mu_2(s,h)=\phi_2(s)\psi_2(h), \quad s\in \mathbb{R}^n, h\in H,
\end{equation}
where $\phi_j(0)=\psi_j(0)=1$, $j=1,2$.

Substituting $s = -\varepsilon s', h'=h$ into (\ref{f1}), we get the
equation

\begin{equation}\label{f4}
\hat\mu_1((I-\varepsilon)s', 2h)\hat\mu_2(0, (I+\varepsilon)h) =
\hat\mu_1(-(I+\varepsilon)s', 0)\hat\mu_2(-2\varepsilon s',
(I-\varepsilon)h), \quad (s,h), \ (s',h') \in Y.
\end{equation}

If $k=1$, i.e. $2h=0$ then equation (\ref{f4}) is of the form

\begin{equation}\label{f5}
\hat\mu_1((I-\varepsilon)s', 0)\hat\mu_2(0, (I+\varepsilon)h) =
\hat\mu_1(-(I+\varepsilon)s', 0)\hat\mu_2(-2\varepsilon s',
(I-\varepsilon)h), \quad (s,h), \ (s',h') \in Y.
\end{equation}
It follows from (\ref{f2}) that $\hat\mu_1(-(I+\varepsilon)s',
0)\neq 0$. Since the condition (\ref{t4}) the equality
$-2\varepsilon (\mathbb{R}^n)=\mathbb{R}^n$ are fulfilled, we obtain
from (\ref{f5}) representation (\ref{f3}) for $\hat\mu_2(s,h)$.

Substituting $s'=-s, h=\varepsilon h'$ into (\ref{f1}), we get the
equation

\begin{equation}\label{f6}
 \hat\mu_1(0,(I+\varepsilon)h)\hat\mu_2((I-\varepsilon)s,2\varepsilon h') =
\hat\mu_1(2s,-(I-\varepsilon)h')\hat\mu_2((I+\varepsilon)s,0), \quad
(s,h), \ (s',h') \in Y.
\end{equation}
It follows from (\ref{f6}) that for $k=1$, i.e. $2h'=0$, the
function $\hat\mu_1(s,h)$ is of the form (\ref{f3}). Thus, the
statement is proved for $k=1$.

Assume that (\ref{f3}) holds if $h$ has order $2^k$. Let $h$ have
order $2^{k+1}$. Then $2h$ has order $2^{k}$, and we have in
(\ref{f4})
$\hat\mu_1((I-\varepsilon)s',2h)=\phi_1((I-\varepsilon)s')\psi_1(2h)$
by induction hypothesis. Arguing similarly to the case of $k=1$, we
obtain from (\ref{f4}) representation (\ref{f3}) for
$\hat\mu_2(s,h)$. Similarly we obtain from (\ref{f6}) the
representation (\ref{f3}) for $\hat\mu_1(s,h)$. The function
$\phi_1(s)$ is a characteristic function of a distribution $\gamma_1
\in \Gamma({\mathbb R}^n)$, and the function $\psi_1(h)$ is a
characteristic function of a distribution $\rho_1$ such that
$\sigma(\rho_1)\subset G$. Thus, $\mu_1=\gamma_1*\rho_1$. Similarly
we obtain that $\mu_2=\gamma_2*\rho_2$. $\blacksquare$

The following corollary follows from the proof of Lemma \ref{t1}.

\begin{corollary}\label{corollary 1}
Let $Y={\mathbb R}^n\times H$, where $H$ is a compact 2-primary
Abelian group. All solutions of equation (\ref{l2}) on $Y$ have form
(\ref{f3}), where the functions $\phi_j(s)$ have form (\ref{f2}).
\end{corollary}

\begin{lemma}\label{compconnected} (\cite{FeTVP1}) Let
$Y$ be a connected compact Abelian group, $X=Y^*$. Let
$\varepsilon\in{\rm Aut}(Y)$ and $(I+\varepsilon)Y=Y$. Let $\mu_1$
and $\mu_2$ be distributions on the group $X$ such that the
characteristic functions $\hat\mu_j(y)$ satisfy equation
$(\ref{l2})$. Then $\hat\mu_j(y)=(x_j, y)$, where $x_j\in X$, $j=1,
2$.
\end{lemma}

\begin{lemma}\label{support}
    Let $X=\mathbb{R}^n\times D$, where $D$ is a countable discrete Abelian group.
    Let $\delta\in Aut(X)$ such that condition $(\ref{i2})$ is fulfilled.
    Let $\xi_1$ and $\xi_2$ be independent random variables with values in $X$ and distributions $\mu_1$ and $\mu_2$.
    If the conditional distribution of the linear statistic $L_2 = \xi_1 + \delta\xi_2$ given $L_1 = \xi_1 + \xi_2$ is symmetric, then
    the random variables $\xi_j$ can be replaced by their shifts $\xi'_j$ with distributions $\mu'_j$ in such a manner that
    $\sigma(\mu'_j) \subset R^n\times b_D$ and the conditional distribution of the linear statistic
    $L'_2 = \xi'_1 + \delta\xi'_2$ given $L'_1 = \xi'_1 + \xi'_2$
    is symmetric.
\end{lemma}

Lemma \ref{support} was actually proved in the article
\cite{FeTVP2}, but was not formulated as a separate statement. For
the sake of exposition, we present its proof in the article.

\bigskip

\textbf{Proof of Lemma \ref{support}.} We shall reduce the proof of
the lemma to the case when the subgroup $D$ is torsion. The group
$Y=X^*$ is topologically isomorphic to the group ${\mathbb
R}^n\times K$, where $K=D^*$. To avoid introducing new notation we
will suppose that $Y={\mathbb R}^n\times K$. Denote by $x=(t, d)$,
where $t\in {\mathbb R}^n$, $d\in D$, elements of the group $X$.
Denote by $y=(s, k)$, where $s\in {\mathbb R}^n$, $k\in K$, elements
of the group $Y$.

It follows from Lemma \ref{l1} that the characteristic functions
$\hat\mu_j(y)$ satisfy equation (\ref{l2}). Since the subgroup $K$
consists of all compact elements of the group $Y$, it is
characteristic. Hence, the subgroup $c_K$ is also characteristic,
i.e. $\varepsilon(c_K)=c_K$. Consider the restriction of equation
(\ref{l2}) on the subgroup $c_K$.

We shall show that
\begin{equation}\label{t2.3}
(I+\varepsilon)(c_K)=c_K.
\end{equation}
Since $A(X,c_K)= \mathbb{R}^n\times b_D$, we get that $(c_K)^*
\approx X/ \mathbb{R}^n\times b_D$. Since the subgroup ${\mathbb
R}^n\times b_D$ is characteristic, the automorphism $\delta$ induces
an automorphism $\bar{\delta}$ on the factor group $X/
\mathbb{R}^n\times b_D$.

The condition (\ref{t2.3}) is equivalent to the condition
\begin{equation}\label{t2.4}
Ker(I+\bar{\delta})=\{0\}.
\end{equation}

If $(t_0, d_0)\in X$ and $[(t_0, d_0)]\in Ker(I+\bar{\delta})$ then
$[(I+\delta)(t_0, d_0)]=[0]$. Hence, $(I+\delta)(t_0, d_0)\in
{\mathbb R}^n\times b_D$. Since the subgroup $b_D$ is torsion, we
get that $k(I+\delta)(t_0, d_0)\in{\mathbb R}^n$ for some natural
$k$. It follows from this that $(I+\delta)k(t_0, d_0)\in{\mathbb
R}^n$, i.e.
\begin{equation}\label{t2.5}
(I+\delta)(kt_0, kd_0)=(t_1, 0)
\end{equation}
for some $t_1\in \mathbb{R}^n$. It is obvious that
$(I+\delta)(\mathbb{R}^n)\subset\mathbb{R}^n$. It follows from
(\ref{i2}) that the restriction of the continuous endomorphism
$I+\delta$ of the group $X$ into the subgroup $\mathbb{R}^n$ is a
topological automorphism of the group $\mathbb{R}^n$. Therefore
\begin{equation}\label{t2.6}
(t_1, 0)=(I+\delta)(t_2, 0)
\end{equation}
for some $t_2\in \mathbb{R}^n$. In view of (\ref{i2}), it follows
from (\ref{t2.5}) and (\ref{t2.6}) that $kd_0=0$. Hence $d_0\in b_D$
and $(t_0, d_0)\in{\mathbb R}^n\times b_D$. It follows from this
that $[(t_0, d_0)]=0$. Thus, (\ref{t2.4}) is fulfilled. So,
(\ref{t2.3})  is fulfilled too.

Lemma \ref{compconnected} implies that $\hat\mu_j(y)=(x_j,y)$, $y\in
c_K$, $j=1, 2$. By the theorem on the extension of a character from
the closed subgroup to the group, we can suppose that $x_j\in X$,
$j=1, 2$. Substituting these expressions of $\hat\mu_j(y)$ into
equation (\ref{l2}) and taking into account the equality $A(X,
c_K)={\mathbb R}^n\times b_D$, we get
\begin{equation}\label{t2.1}
2(x_1+\delta x_2)\in \mathbb{R}^n\times b_D.
\end{equation}
Since $b_D$ consists of elements of finite order of the group $X$,
it follows from (\ref{t2.1}) that
\begin{equation}\label{t2.2}
x_1+\delta x_2\in \mathbb{R}^n\times b_D.
\end{equation}

Consider new random variables $\xi'_1=\xi_1+\delta x_2$ и
$\xi'_2=\xi_2-x_2$ with values in the group $X$. Denote by $\mu'_j$
distributions of the random variables $\xi'_j$. Then
$\mu'_1=\mu_1*E_{\delta x_2}$, $\mu'_2=\mu_2*E_{-x_2}$. It is easy
to see that the characteristic functions $\hat\mu'_j(y)$ satisfy
equation (\ref{l2}). Then Lemma \ref{l1} implies that the
conditional distribution of the linear statistic $L'_2 = \xi'_1 +
\delta\xi'_2$ given $L'_1 = \xi'_1 + \xi'_2$ is symmetric. We have
that $\hat\mu'_2(y)=1$, $y\in c_K$. In view of (\ref{t2.2}), we get
that $\hat\mu'_1(y)=1$, $y\in c_K$. It follows from Lemma \ref{supp}
that $\sigma(\mu'_j)\subset A(X, c_K)=\mathbb{R}^n\times b_D$, $j=1,
2$. $\blacksquare$

\begin{lemma}\label{fel2} (\cite{FeTVP2})
Let $X$ be a countable discrete Abelian group, $G$ be a subgroup
generated by all elements of odd order of the group $X$. Let
$\delta\in Aut(X)$ such that condition $(\ref{i2})$ is fulfilled.
Let $\xi_1$ and $\xi_2$ be independent random variables with values
in $X$ and distributions $\mu_1$ and $\mu_2$. If the conditional
distribution of the linear statistic $L_2 = \xi_1 + \delta\xi_2$
given $L_1 = \xi_1 + \xi_2$ is symmetric, then
$\mu_j=\rho_j*m_K*E_{x_j}$, where $\sigma(\rho_j)\subset X_2$,  $K$
is a finite subgroup of $G$, $x_j\in X$, $j=1, 2$.
\end{lemma}

\begin{lemma}\label{fel3} (\cite{FeTVP2}) Let $X$ be a countable discrete Abelian group without elements of order $2$.
Let $\delta\in Aut(X)$ such that condition $(\ref{i2})$ is
fulfilled. Let $\xi_1$ and $\xi_2$ be independent random variables
with values in $X$ and distributions $\mu_1$ and $\mu_2$. If the
conditional distribution of the linear statistic $L_2 = \xi_1 +
\delta\xi_2$ given $L_1 = \xi_1 + \xi_2$ is symmetric, then
$\mu_j=m_F*E_{x_j}$, where $F$ is a finite subgroup of $X$, $x_j\in
X$, $j=1, 2$. Moreover, $\delta (F)=F$.
\end{lemma}

\textbf{Proof of Theorem \ref{m}.} Since the subgroup ${\mathbb
R}^n\times b_D$ of the group $X$ is characteristic, Lemma
\ref{support} implies that it suffices to prove theorem in the case
where  $D$ is a  torsion subgroup.

A discrete torsion group can be decomposed into a weak direct
product of its $p$-primary components:
$D=P\hspace{-2.2ex}\raisebox{-1.7 ex}{$\scriptscriptstyle p \in
\mathcal{P} $} \ D_p, $ where $\mathcal{P}$ is the set of prime
numbers (see \cite[Theorem 8.4]{Fu1}). Put $G=D_2, \
L=\sum\hspace{-2.2ex}\raisebox{-1.7 ex}{$\scriptscriptstyle p>2 $} \
D_p$. Then $X=\mathbb{R}^n\times G\times L$ and $Y\approx
\mathbb{R}^n\times H\times M$, where $H=G^*, M=L^{*}$. To avoid
introducing new notation we will assume that $Y=\mathbb{R}^n\times
H\times M$. Denote by $(s,h,m), \ s \in \mathbb{R}^n, h \in H, m \in
M,$ elements of the group $Y$. Since $\mathbb{R}^n, H, M$ are
characteristic subgroups, each automorphism $\varepsilon \in {\rm
Aut(Y)}$ can be written in the form $\varepsilon(s,h,m)=
(\varepsilon s, \varepsilon h, \varepsilon m), \ (s,h,m) \in Y$.

It follows from Lemma \ref{l1} that the characteristic functions
$\hat\mu_j(y)$ satisfy equation (\ref{l2}). Consider the restriction
of equation (\ref{l2}) to the subgroup $M$.  Lemma \ref{fel3} and
(\ref{11a}) imply that

\begin{equation}\label{t2.7}
\hat\mu_1(0,0,m)=\left\{%
\begin{array}{ll}
    (l_1,m), & \hbox{$m\in A(M,F)$;} \\
    0, & \hbox{$m\not\in A(M,F)$;} \\
\end{array}%
\right.
\hat\mu_2(0,0,m)=\left\{%
\begin{array}{ll}
    (l_2,m), & \hbox{$m\in A(M,F)$;} \\
    0, & \hbox{$m\not\in A(M,F)$;} \\
\end{array}%
\right.
\end{equation}
\noindent where $l_j\in L$, $F$ is a finite subgroup of $L$ such
that $\delta(F)=F$. Substituting (\ref{t2.7}) into equation
(\ref{l2}) and considering  the restriction of equation (\ref{l2})
to the subgroup $M$, we get
\begin{equation}\label{t2.8}
2(l_1+\delta l_2)\in F.
\end{equation}
Since $L$ does not contain elements of order 2, it follows from
(\ref{t2.8}) that
\begin{equation}\label{t2.9}
l_1+\delta l_2\in F.
\end{equation}
Considering new random variables $\zeta_1=\xi_1+\delta l_2$ and
$\zeta_2=\xi_2-l_2$ and reasoning as in the end of the proof of
Lemma \ref{support}, we obtain that we can suppose from the
beginning that

\begin{equation}\label{t2.10}
\hat\mu_1(0,0,m)=\left\{%
\begin{array}{ll}
    1, & \hbox{$m\in A(M,F)$;} \\
    0, & \hbox{$m\not\in A(M,F)$;} \\
\end{array}%
\right.
\hat\mu_2(0,0,m)=\left\{%
\begin{array}{ll}
    1, & \hbox{$m\in A(M,F)$;} \\
    0, & \hbox{$m\not\in A(M,F)$.} \\
\end{array}%
\right.
\end{equation}
It follows from this that $\sigma(\mu_j)\subset A(X,
A(M,F))=\mathbb{R}^n\times G\times F$, $j=1, 2$.

Since the subgroups $\mathbb{R}^n$ and $G$ are characteristic and
$\delta(F)=F$, it suffices to prove the theorem in the case when
$X=\mathbb{R}^n\times G\times L$, where the subgroup $L$ is finite
and $\delta(L)=L$. Moreover, the following representations are valid

\begin{equation}\label{t2.19}
\hat\mu_1(0,0,m)=\left\{%
\begin{array}{ll}
    1, & \hbox{$m= 0$;} \\
    0, & \hbox{$m\neq 0$;} \\
\end{array}%
\right.
\hat\mu_2(0,0,m)=\left\{%
\begin{array}{ll}
    1, & \hbox{$m=0$;} \\
    0, & \hbox{$m\neq 0$.} \\
\end{array}%
\right.
\end{equation}

Putting $u=v=(0,0,m)$ into (\ref{l2}) $u=v=(0,0,m)$, we get

\begin{equation}\label{t2.20}
\hat\mu_1(0, 0, 2m)\hat\mu_2(0, 0, (I+\varepsilon) m) =\hat\mu_2(0,
0, (I-\varepsilon) m), \quad   m \in M.
\end{equation}

Since the subgroup $M$ does not contain elements of order 2, $2m=0$
if and only if $m=0$. Then $\hat\mu_1(0, 0, 2m)=0$ for $m\neq 0$,
and it follows from (\ref{t2.20}) that $\hat\mu_2(0, 0,
(I-\varepsilon) m)=0$ if and only if $m\neq 0$. It follows from this
and the representations (\ref{t2.19}) that the restriction of the
endomorphism $I-\varepsilon$ into $M$ has a zero kernel. Since the
subgroup $M$ is finite, the restriction of the endomorphism
$I-\varepsilon$ into $M$ is an automorphism of $M$, i.e.

\begin{equation}\label{t2.21}
(I-\varepsilon) M=M.
\end{equation}

Consider the restriction of equation (\ref{l2}) on the subgruop
$H\times M$. Lemma \ref{fel2} implies that

\begin{equation}\label{t2.11}
\hat\mu_1(0,h,m)=\left\{%
\begin{array}{ll}
    \psi_1(h), & \hbox{$m=0$;} \\
    0, & \hbox{$m\neq 0$;} \\
\end{array}%
\right.
\hat\mu_2(0,h,m)=\left\{%
\begin{array}{ll}
    \psi_2(h), & \hbox{$m=0$;} \\
    0, & \hbox{$m\neq 0$.} \\
\end{array}%
\right.
\end{equation}
\noindent where $\psi_j(h)$ are characteristic functions on $H$.

Rewrite equation (\ref{l2}) in the form
\begin{equation}\label{t2.12}
\hat\mu_1(s+s', h+h', l+l')\hat\mu_2(s+\varepsilon s', h+\varepsilon
h', m+\varepsilon m') =$$ $$ = \hat\mu_1(s-s', h-h',
l-l')\hat\mu_2(s-\varepsilon s', h-\varepsilon h', m-\varepsilon
m'), \quad\quad (s,h,m), \ (s',h',m') \in Y.
\end{equation}

Putting $s'=s,h'=-h,m'=-m$ into (\ref{t2.12}), we get

\begin{equation}\label{t2.13}
\hat\mu_1(2s,0,0)\hat\mu_2((I+\varepsilon)s,(I-\varepsilon)h,(I-\varepsilon)m)
=$$$$=
\hat\mu_1(0,2h,2m)\hat\mu_1((I-\varepsilon)s,(I+\varepsilon)h,(I+\varepsilon)m),
\quad\quad (s,h,m) \in Y.
\end{equation}
It follows from (\ref{t2.11}) that $\hat\mu_1(0,2h,2m)=0$ for $m\neq
0$. Hence,
$\hat\mu_1(2s,0,0)\hat\mu_2((I+\varepsilon)s,(I-\varepsilon)h,(I-\varepsilon)m)=0$
for $m\neq 0$. Consider the restriction of equation (\ref{l2}) on
the subgroup $\mathbb{R}^n$. It was proved in \cite{FeTVP1} that all
solutions of this equation are the characteristic functions of
Gaussian distributions, i.e. have form (\ref{f2}). It follows from
(\ref{f2}) that the function $\hat\mu_1(2s,0,0)$ do not vanish. Then
we get from (\ref{t2.13}) that

\begin{equation}\label{t2.14}
\hat\mu_2((I+\varepsilon)s,(I-\varepsilon)h,(I-\varepsilon)m)=0,
\quad s\in \mathbb{R}^n, h\in H,
\end{equation}

\noindent for $m\neq 0$. Note that it follows from $(\ref{i2})$ that
$(I+\varepsilon)\mathbb{R}^n=\mathbb{R}^n$. As in the proof of the
part 1 of Lemma \ref{t1}, we obtain (\ref{t4}). Taking into account
(\ref{t2.21}), it follows from (\ref{t2.14}) that

\begin{equation}\label{t2.15}
\hat\mu_2(s,h,m)=0, \quad s\in \mathbb{R}^n, h\in H, m\neq 0.
\end{equation}

\noindent Similarly we get that

\begin{equation}\label{t2.16}
\hat\mu_1(s,h,m)=0, \quad s\in \mathbb{R}^n, h\in H, m\neq 0.
\end{equation}

Put $m=m'=0$ in (\ref{t2.12}). Corollary \ref{corollary 1} implies
that

\begin{equation}\label{t2.17}
 \hat\mu_1(s,h,0) = \phi_1(s)\psi_1(h), \quad \hat\mu_2(s,h,0)
= \phi_2(s)\psi_2(h),
\end{equation}
where the functions $\phi_1(s),\phi_2(s)$ are of the form
(\ref{f2}), and the functions $\psi_1(h),\psi_2(h)$ are the
characteristic functions of distributions $\rho_j$ such that
$\sigma(\rho_j)\subset G$.

We deduce from (\ref{t2.15})-(\ref{t2.17}) that

\begin{equation}\label{t2.18}
\hat\mu_1(s,h,m)=\left\{%
\begin{array}{ll}
    \phi_1(s)\psi_1(h), & \hbox{$m=0$;} \\
    0, & \hbox{$m\neq 0$;} \\
\end{array}%
\right.
\hat\mu_2(s,h,m)=\left\{%
\begin{array}{ll}
    \phi_2(s)\psi_2(h), & \hbox{$m=0$;} \\
    0, & \hbox{$m\neq 0$.} \\
\end{array}%
\right.
\end{equation}
It is easy to see that the assertion of the theorem follows from
(\ref{t2.18}). $\blacksquare$


\end{document}